\documentclass[letterpaper, 10 pt, conference]{IEEEtran}  

\IEEEoverridecommandlockouts                              

\usepackage{epsfig} 
\usepackage{epstopdf}
\usepackage{amsmath} 
\usepackage{amssymb}  
\usepackage{caption}

\usepackage[caption=false,font=footnotesize]{subfig}

\DeclareCaptionLabelSeparator{periodspace}{.\quad}
\title{\LARGE \bf
Attitude and Angular Velocity Tracking for a Rigid Body using Geometric Methods on the Two-Sphere (Stability Proof)
}

\author{Michalis Ramp$^{1}$ and Evangelos Papadopoulos$^{2}$ 
\thanks{*Financial support by the European Union (European Social Fund-ESF) and Greek national funds through the Operational Program ''Education and Lifelong Learning'' of  the National Strategic Reference Framework Research Funding Program: THALES: Reinforcement of the interdisciplinary and/or interinstitutional research $\&$ innovation is acknowledged.}
\thanks{$^{1}$M. Ramp is with the Department of Mechanical Engineering, National Technical University of Athens, (NTUA) 15780 Athens, Greece
        {\tt\small rampmich@mail.ntua.gr}}%
\thanks{$^{2}$E. Papadopoulos is with the Department of Mechanical Engineering, NTUA, 15780 Athens (tel: +30-210-772-1440; fax: +30-210-772-1455)
        {\tt\small egpapado@central.ntua.gr}}%
}

\begin{document}


\IEEEpubid{\copyright~2015 IEEE. Personal use is permitted. For any other purposes, permission must be obtained from the IEEE. DOI: 10.1109/ECC.2015.7331033}
\IEEEpubidadjcol
\maketitle
\pagestyle{empty}

\begin{abstract}
Stability proof of the controller proposed in \cite{Ramp}.
This proof must be studied together with Sec. III-C in \cite{Ramp}.
\end{abstract}


\section{Pointing direction and angular velocity tracking}
A control system is defined in $\text{L}_{2}$ under the assumption that we have  a fairly accurate estimate of the model parameters showing exponential convergence in an envelope around the zero equilibrium of ${}^{b}\mathbf{e}_{q}$, ${}^{b}\mathbf{e}_{\omega}$ using Lyapunov analysis.

\textit{Proposition 1:} For $\Lambda,\gamma,\eta\in\mathbb{R}^{+}$ and a desired pointing direction curve $\mathbf{q}_{d}(t)\in\text{S}^{2}$ and a desired angular velocity profile $^{b}\boldsymbol{\omega}_{d}(t)$, $^{b}\dot{\boldsymbol{\omega}}_{d}(t)$, around the $\mathbf{q}_{d}(t)$ axis, we define the control moment $^{b}\mathbf{u}$ as,
\begin{IEEEeqnarray}{rCl}
\label{eq:u}
^{b}\mathbf{u}&{=}&\eta^{-1}\widehat{\mathbf{J}}(-\eta(\widehat{\mathbf{f}}+\mathbf{d}){-}(\Lambda+\Psi){}^{b}\dot{\mathbf{e}}_{q}{-}\dot{\Psi}{}^{b}{\mathbf{e}}_{q}{-}\gamma\mathbf{s})\IEEEyessubnumber\label{eq:uc}\\
\mathbf{d}&{=}&({}^{b}\boldsymbol{\omega})^{\times}\mathbf{Q}^{T}\mathbf{Q}_{d}{}^{b}\boldsymbol{\omega}_{d}-\mathbf{Q}^{T}\mathbf{Q}_{d}{}^{b}\dot{\boldsymbol{\omega}}_{d}\IEEEyessubnumber\label{eq:d}\\
\mathbf{f}&{=}&\mathbf{J}^{-1}((\mathbf{J}{}^{b}\boldsymbol{\omega})^{\times}({}^{b}\boldsymbol{\omega})-c{}^{b}\boldsymbol{\omega}+\boldsymbol{\tau})\IEEEyessubnumber\\
\mathbf{s}&{=}&(\Lambda+\Psi){}^{b}\mathbf{e}_{q}+\eta{}^b\mathbf{e}_{\boldsymbol{\omega}}\IEEEyessubnumber\label{eq:s}\\
\gamma&=&\gamma_{1}+\gamma_{2}+\gamma_{3},\gamma_{3}=\gamma_{4}+\gamma_{5},\gamma_{i=1,..,5}\in\mathbb{R}^{+}\IEEEyessubnumber\label{eq:cond0}\\
\gamma_{3}&<&\frac{\gamma_{4}(\Lambda+\Psi)^{2}}{\Psi^{2}}\IEEEyessubnumber\label{eq:cond0.5}\\
\gamma_{5}&>&A_{2_{max}}/\lambda_{J}\IEEEyessubnumber\label{eq:cond1}\\
\eta&>&A_{1_{max}}/(\gamma_{5}\lambda_{J}-A_{2_{max}})\IEEEyessubnumber\label{eq:cond2}
\end{IEEEeqnarray}
where  $\lambda_J=\lambda_{min}(\mathbf{J}^{-1}\widehat{\mathbf{J}})$ while $\widehat{(.)}$ signifies estimated parameters due to parameter identification errors. 
The bounds $A_{1_{max}},A_{2_{max}}$ are found via the expressions after eq. \eqref{eq:dVq} and are given in \eqref{eq:bounds}.
It will be shown that the above control law stabilizes and maintains ${}^{b}\mathbf{e}_{q}$, ${}^{b}\mathbf{e}_{\omega}$ in a bounded set around the zero equilibrium.
Furthermore for perfect knowledge of the system parameters the above law stabilizes  ${}^{b}\mathbf{e}_{q}$, ${}^{b}\mathbf{e}_{\omega}$ to zero exponentially.

\textit{Proof:} We utilize a sliding structure in $\text{L}_{2}$ by defining the surface in terms of the configuration error vectors (12), (15) and the attitude error function (8) so that they appear explicitly in the Lyapunov candidate function.
Then the control design is similar to nonlinear control design in Euclidean spaces  \cite{GCBOOK}, \cite{KHALIL}.
The defined sliding surface is given in \eqref{eq:s}. Its derivative is,
\begin{IEEEeqnarray*}{l}
\label{eq:ds}
\dot{\mathbf{s}}=\dot{\Psi}{}^{b}\mathbf{e}_{q}+(\Lambda+\Psi){}^{b}\dot{\mathbf{e}}_{q}+\eta{}^{b}\dot{\mathbf{e}}_{\omega}\IEEEyesnumber
\end{IEEEeqnarray*}
For $\kappa=2\eta\Lambda(\gamma_{2}+\gamma_{3})\lambda_J$ the Lyapunov candidate is,
\begin{IEEEeqnarray}{C}
\label{eq:V}
V(\Psi,{}^{b}\mathbf{e}_{q},{}^{b}\mathbf{e}_{\omega})=\frac{1}{2}\mathbf{s}^{T}\mathbf{s}+\kappa\Psi
\end{IEEEeqnarray}
Using (13) and the vector $\mathbf{z}_{q}=[\lVert{}^{b}{\mathbf{e}}_{q}\rVert;\lVert{}^{b}{\mathbf{e}}_{\omega}\rVert]$ it holds that,
\begin{IEEEeqnarray*}{C}
\lambda_{min}(\mathbf{W}_{1})\lVert\mathbf{z}_{q}\rVert^{2}{\leq}\mathbf{z}_{q}^{T}\mathbf{W}_{1}\mathbf{z}_{q}{\leq}{V}{\leq}\mathbf{z}_{q}^{T}\mathbf{W}_{2}\mathbf{z}_{q}{\leq}\lambda_{max}(\mathbf{W}_{2})\lVert\mathbf{z}_{q}\rVert^{2}\\
\mathbf{W}_{1}=\begin{bmatrix}
\frac{(\Lambda+\Psi)^{2}}{2}+\kappa&-\frac{(\Lambda+\Psi)\eta}{2}\\
-\frac{(\Lambda+\Psi)\eta}{2}&{\eta^{2}}/{2}
\end{bmatrix}\\
\mathbf{W}_{2}=\begin{bmatrix}
\frac{(\Lambda+\Psi)^{2}}{2}+2\kappa&\frac{(\Lambda+\Psi)\eta}{2}\\
\frac{(\Lambda+\Psi)\eta}{2}&{\eta^{2}}/{2}
\end{bmatrix}
\end{IEEEeqnarray*}
Differentiating (\ref{eq:V}) and substituting (\ref{eq:ds}) we get,
\begin{IEEEeqnarray}{rCl}
\dot{V}&{=}&\mathbf{s}^{T}\dot{\mathbf{s}}+\kappa\dot{\Psi}\IEEEnonumber\\
&{=}&\mathbf{s}^{T}(({}^{b}\mathbf{e}^T_{q}{}^{b}\mathbf{e}_{\omega}){}^{b}\mathbf{e}_{q}{+}(\Lambda{+}\Psi){}^{b}\dot{\mathbf{e}}_{q}{+}\eta\mathbf{J}^{-1}{}^{b}\mathbf{u}{+}\eta\mathbf{f}{+}\eta\mathbf{d})\label{eq:dV}\\
&&+\kappa{}^{b}\mathbf{e}^T_{q}{}^{b}\mathbf{e}_{\omega}\IEEEnonumber
\end{IEEEeqnarray}
To avoid high frequency chattering and the discontinuities introduced by the standard sliding condition, and since a control system is developed on $\text{S}^{2}$, the conventional sliding condition will not be used.
Instead, the control law is designed such that when not on the surface, the following holds,
\begin{IEEEeqnarray*}{C}
\mathbf{s}^{T}\dot{\mathbf{s}}\leqq-k\lVert\mathbf{s}\rVert^{2},\;k>0
\end{IEEEeqnarray*}
Substituting \eqref{eq:uc} and \eqref{eq:deq} to (\ref{eq:dV}), after some manipulations,
\begin{IEEEeqnarray*}{rCl}
\dot{V}&{=}&\mathbf{s}^{T}(\Delta\mathbf{J}\left(({}^{b}\mathbf{e}_{q}{}^{b}\mathbf{e}^T_{q}){}^{b}\mathbf{e}_{\omega}{+}(\Lambda{+}\Psi)(\mathbf{E}{}^{b}\mathbf{e}_{\omega}{+}\boldsymbol{\Xi}{}^{b}\boldsymbol{\omega}_{d})\right)\\
&&{+}\eta(\mathbf{f}{-}\mathbf{J}^{-1}\widehat{\mathbf{J}}\widehat{\mathbf{f}})+\Delta\mathbf{J}\eta\mathbf{d}{-}\mathbf{J}^{-1}\widehat{\mathbf{J}}\gamma\mathbf{s})+\kappa{}^{b}\mathbf{e}^T_{q}{}^{b}\mathbf{e}_{\omega}\\
\Delta\mathbf{J}&=&\mathbf{I}-\mathbf{J}^{-1}\widehat{\mathbf{J}}
\end{IEEEeqnarray*}
Employing \eqref{eq:dalt}, after several manipulations,
\begin{IEEEeqnarray}{rCl}
\label{eq:dVq}
\dot{V}&=&\mathbf{s}^{T}(-\gamma\mathbf{J}^{-1}\widehat{\mathbf{J}}\mathbf{s}+\mathbf{A}{}^{b}\mathbf{e}_{\omega}+\mathbf{B}){+}\kappa{}^{b}\mathbf{e}_{q}^{T}{}^{b}\mathbf{e}_{\omega}
\end{IEEEeqnarray}
where $\mathbf{A}\in\mathbb{R}^{3{\times}3}$, $\mathbf{B}\in\mathbb{R}^{3{\times}1}$, are given by,
\begin{IEEEeqnarray*}{rCl}
\mathbf{A}&=&\mathbf{A}_{1}-\eta\mathbf{A}_{2}\\
\mathbf{A}_{1}&=&\Delta\mathbf{J}\left({}^{b}\mathbf{e}_{q}{}^{b}\mathbf{e}^T_{q}{+}(\Lambda{+}\Psi)\mathbf{E}\right),
\mathbf{A}_{2}=\Delta\mathbf{J}(\mathbf{Q}^{T}\mathbf{Q}_{d}{}^{b}\boldsymbol{\omega}_{d})^{\times}\\
\mathbf{B}&=&\Delta\mathbf{J}\left((\Lambda{+}\Psi)\boldsymbol{\Xi}{}^{b}\boldsymbol{\omega}_{d}-\eta\mathbf{Q}^{T}\mathbf{Q}_{d}{}^{b}\dot{\boldsymbol{\omega}}_{d}\right)+\eta(\mathbf{f}{-}\mathbf{J}^{-1}\widehat{\mathbf{J}}\widehat{\mathbf{f}})
\end{IEEEeqnarray*}
\IEEEpubidadjcol
Under the assumption mentioned in the beginning of this section, i.e. we have  a fairly accurate estimate of the model parameters, the following holds,
\begin{IEEEeqnarray*}{rCl}
\lVert\mathbf{f}{-}\mathbf{J}^{-1}\widehat{\mathbf{J}}\widehat{\mathbf{f}}\rVert\leq f_{max}<\infty,f_{max}\in\mathbb{R}^{+}
\end{IEEEeqnarray*}
Additionally since ${}^{b}{\boldsymbol{\omega}}_{d},{}^{b}\dot{\boldsymbol{\omega}}_{d}$ are bounded, the following holds,
\begin{IEEEeqnarray}{C}
\exists A_{1_{max}},A_{2_{max}},B_{_{max}}\in\mathbb{R}^{+}-\{\infty\}\IEEEnonumber\\
\lVert\mathbf{A}_{1}\rVert\leq A_{1_{max}},\lVert\mathbf{A}_{2}\rVert\leq A_{2_{max}},\lVert\mathbf{B}\rVert\leq B_{_{max}}\label{eq:bounds}
\end{IEEEeqnarray}
Expanding (\ref{eq:dVq}) and rearranging we have,
\begin{IEEEeqnarray*}{rCl}
\dot{V}&\leq&-\gamma\lambda_{J}\mathbf{s}^{T}\mathbf{s}+\Upsilon+\breve{\mathbf{A}}\lVert{}^{b}\mathbf{e}_{\omega}\rVert+(\eta^{2}\lVert\mathbf{A}_{2}\rVert+\eta\lVert\mathbf{A}_{1}\rVert)\lVert{}^{b}\mathbf{e}_{\omega}\rVert^{2}\\
&&{+}\kappa{}^{b}\mathbf{e}_{q}^{T}{}^{b}\mathbf{e}_{\omega}\\
&\leq&-\gamma\lambda_{J}\mathbf{s}^{T}\mathbf{s}+\Upsilon+\breve{\mathbf{A}}\lVert{}^{b}\mathbf{e}_{\omega}\rVert+(\eta^{2}A_{2_{max}}+\eta A_{1_{max}})\lVert{}^{b}\mathbf{e}_{\omega}\rVert^{2}\\
&&{+}\kappa{}^{b}\mathbf{e}_{q}^{T}{}^{b}\mathbf{e}_{\omega}\\
\Upsilon&=&(\Lambda{+}\Psi)\lVert\mathbf{B}\rVert\\
\breve{\mathbf{A}}&=&\eta\lVert\mathbf{B}\rVert+(\Lambda+\Psi)\lVert\mathbf{A}\rVert
\end{IEEEeqnarray*}
Using \eqref{eq:cond0} after several manipulations,
\begin{IEEEeqnarray*}{rCl}
\dot{V}&\leq&-\gamma_{1}\lambda_{J}\mathbf{s}^{T}\mathbf{s}-\mathbf{z}_{q}^{T}\mathbf{W}_{3}\mathbf{z}_{q}+\Upsilon-\gamma_{3}\lambda_{J}(\Lambda+\Psi)^{2}\lVert{}^{b}\mathbf{e}_{q}\rVert^{2}\\
&&-\gamma_{3}\lambda_{J}2\Psi\eta{}^{b}\mathbf{e}^{T}_{q}{}^{b}\mathbf{e}_{\omega}-\gamma_{3}\lambda_{J}\eta^{2}\lVert{}^{b}\mathbf{e}_{\omega}\rVert^{2}\\
&&+\breve{\mathbf{A}}\lVert{}^{b}\mathbf{e}_{\omega}\rVert+(\eta^{2}A_{2_{max}}+\eta A_{1_{max}})\lVert{}^{b}\mathbf{e}_{\omega}\rVert^{2}\\
&&\mathbf{W}_{3}=\begin{bmatrix}
\gamma_{2}\lambda_{J}(\Lambda+\Psi)^{2}&-\gamma_{2}\lambda_{J}\Psi\eta\\
-\gamma_{2}\lambda_{J}\Psi\eta&\eta^{2}\gamma_{2}\lambda_{J}
\end{bmatrix}
\end{IEEEeqnarray*}
Employing $\gamma_{3}=\sum_{i=4}^{5}\gamma_{i},\gamma_{i},\in\mathbb{R}^{+}$ after some manipulations,
\begin{IEEEeqnarray*}{rCl}
\dot{V}&\leq&-\gamma_{1}\lambda_{J}\mathbf{s}^{T}\mathbf{s}-\mathbf{z}_{q}^{T}\mathbf{W}_{3}\mathbf{z}_{q}-\mathbf{z}_{q}^{T}\mathbf{W}_{4}\mathbf{z}_{q}+\Upsilon\\
&&-\left((\gamma_{5}\lambda_{J}- A_{2_{max}})\eta^{2}-\eta A_{1_{max}}\right)\lVert{}^{b}\mathbf{e}_{\omega}\rVert^{2}+\breve{\mathbf{A}}\lVert{}^{b}\mathbf{e}_{\omega}\rVert\\
&&\mathbf{W}_{4}=\begin{bmatrix}
\gamma_{3}\lambda_{J}(\Lambda+\Psi)^{2}&-\gamma_{3}\lambda_{J}\Psi\eta\\
-\gamma_{3}\lambda_{J}\Psi\eta&\gamma_{4}\lambda_{J}\eta^{2}
\end{bmatrix}
\end{IEEEeqnarray*}
Via \eqref{eq:cond0.5} the matrix $\mathbf{W}_{4}$ is positive definite and via the conditions \eqref{eq:cond1}, \eqref{eq:cond2}, the fifth term of the inequality above is non-positive.
For ${}^{b}\mathbf{e}_{\omega}$ such that,
\begin{IEEEeqnarray}{rCl}
\lVert{}^{b}\mathbf{e}_{\omega}\rVert>\frac{\breve{\mathbf{A}}}{(\gamma_{5}\lambda_{J}-A_{2_{max}})\eta^{2}-\eta A_{1_{max}}}\label{eq:beobound}
\end{IEEEeqnarray}
then
\begin{IEEEeqnarray*}{rCl}
\dot{V}&\leq&-\gamma_{1}\lambda_{J}\mathbf{s}^{T}\mathbf{s}-\lambda_{min}(\mathbf{W}_{3})\lVert\mathbf{z}_{q}\rVert^{2}-\lambda_{min}(\mathbf{W}_{4})\lVert\mathbf{z}_{q}\rVert^{2}+\Upsilon
\end{IEEEeqnarray*}
Finally for,
\begin{IEEEeqnarray}{rCl}
\lVert\mathbf{z}_{q}\rVert>\sqrt{\frac{\Upsilon}{\lambda_{min}(\mathbf{W}_{3})}}\label{eq:zqbound}
\end{IEEEeqnarray}
then
\begin{IEEEeqnarray*}{rCl}
\dot{V}&\leq&-\gamma_{1}\lambda_{J}\mathbf{s}^{T}\mathbf{s}-\lambda_{min}(\mathbf{W}_{4})\lVert\mathbf{z}_{q}\rVert^{2}\leq-\lambda_{min}(\mathbf{W}_{4})\lVert\mathbf{z}_{q}\rVert^{2}
\end{IEEEeqnarray*}
Notice that the first term above ensures sliding behavior.

\textit{Boundedness:} Employing conditions \eqref{eq:beobound} and \eqref{eq:zqbound} the following sets are defined,
\begin{IEEEeqnarray*}{rCl}
M_{1}=\{({}^{b}\mathbf{e}_{q},{}^{b}\mathbf{e}_{\omega})\in\mathbb{R}^{3}\times\mathbb{R}^{3}|\text{Eq.} \eqref{eq:zqbound}, \text{Eq.} \eqref{eq:beobound}\}
\end{IEEEeqnarray*}
To ensure that $\lVert{}^{b}\mathbf{e}_{q}\rVert<1$, i.e. the states are in $\text{L}_{2}$, the following set is defined,
\begin{IEEEeqnarray*}{rCl}
M_{2}=\{({}^{b}\mathbf{e}_{q},{}^{b}\mathbf{e}_{\omega})\in\mathbb{R}^{3}\times\mathbb{R}^{3}|\lVert\mathbf{z}_{q}\rVert<1, \text{Eq.} \eqref{eq:beobound}\}
\end{IEEEeqnarray*}
Finally for proper $\gamma,{}^{b}\boldsymbol{\omega}_{d_{max}},{}^{b}\dot{\boldsymbol{\omega}}_{d_{max}}$ the following hold,
\begin{IEEEeqnarray}{rCl}
\Upsilon_{max}&{<}&\lambda_{min}(\mathbf{W}_{3})\IEEEyessubnumber\label{eq:setcond1}\\
\breve{\mathbf{A}}&{<}&\left({(\gamma_{5}\lambda_{J}{-}A_{2_{max}})\eta^{2}{-}\eta A_{1_{max}}}\right)\!\!\sqrt{\!\frac{\Upsilon_{max}}{\lambda_{min}(\mathbf{W}_{3})}}\IEEEyessubnumber\label{eq:setcond2}
\end{IEEEeqnarray}
Conditions \eqref{eq:setcond1}, \eqref{eq:setcond2}, ensure that $M_{1}\subset M_{2}$. 

Thus ${}^{b}\mathbf{e}_{q}$, ${}^{b}\mathbf{e}_{\omega}$ are stabilized exponentially in an envelope of radius $\sqrt{\frac{\Upsilon_{max}}{\lambda_{min}(\mathbf{W}_{3})}}$ around the zero equilibrium, with the radius decreasing as $\gamma$ increases.
Additionally for perfect knowledge of the system parameters, we have perfect cancellation, which yields,
\begin{IEEEeqnarray*}{rCl}
\label{eq:27}
\dot{V}&{\leq}&-\gamma_{1}\mathbf{s}^{T}\mathbf{s}{-}\lambda_{min}(\mathbf{W}_{5})\lVert\mathbf{z}_{q}\rVert^{2}{\leq}{-}\frac{\lambda_{min}(\mathbf{W}_{5})}{\lambda_{max}(\mathbf{W}_{2})}V\IEEEyesnumber\\
&&\mathbf{W}_{5}=\begin{bmatrix}
(\gamma_{2}{+}\gamma_{3})(\Lambda{+}\Psi)^{2}&-(\gamma_{2}{+}\gamma_{3})\Psi\eta\\
-(\gamma_{2}{+}\gamma_{3})\Psi\eta&(\gamma_{2}{+}\gamma_{3})\eta^{2}
\end{bmatrix}
\end{IEEEeqnarray*}
and by the comparison lemma, \cite{KHALIL},
\begin{IEEEeqnarray*}{rCl}
V(t)&{\leqq}&V(0)e^{-\frac{\lambda_{min}(\mathbf{W}_{5})}{\lambda_{max}(\mathbf{W}_{2})}t}
\end{IEEEeqnarray*}
Proving that the zero equilibrium of the attitude and angular velocity tracking errors ${}^{b}\mathbf{e}_{q}$, ${}^{b}\mathbf{e}_{\omega}$, is exponentially stable. $\blacksquare$




\renewcommand{\theequation}{A\arabic{equation}}
\setcounter{equation}{0}  
\section*{APPENDIX}  
Vector space isomorphism where $\mathbf{r}\in\mathbb{R}^3$,
\begin{IEEEeqnarray}{C}
(\mathbf{r})^{\times}{=}[0,-r_3,r_2;r_3,0,-r_1;-r_2,r_1,0],((\mathbf{r})^{\times})^{\vee}{=}\mathbf{r}\IEEEeqnarraynumspace\label{r}
\end{IEEEeqnarray}
Exponential map using the Rodrigues formulation \cite{O'Reilly},
\begin{IEEEeqnarray}{C}
\text{exp}(\epsilon\boldsymbol{\xi}^{\times})=\mathbf{I}+\boldsymbol{\xi}^{\times}\sin{\epsilon}+(\boldsymbol{\xi}^{\times})^{2}(1-\cos{\epsilon})\label{expm}
\end{IEEEeqnarray}
Derivative of the configuration error vector (12),
\begin{IEEEeqnarray}{rCl}
\label{eq:deq}
{}^{b}\dot{\mathbf{e}}_{q}&{=}&\mathbf{E}{}^{b}\mathbf{e}_{\omega}+\boldsymbol{\Xi}{}^{b}\boldsymbol{\omega}_{d}\\
\mathbf{E}&{=}&\frac{({}^{b}\mathbf{e}_{q}\mathbf{q}^{T}_{d})(\mathbf{q})^{\times}\mathbf{Q}}{2(1{+}\mathbf{q}^{T}\mathbf{q}_{d})}+\frac{\mathbf{Q}^{T}\{\mathbf{q}_{d}^{T}\mathbf{q}\mathbf{I}-\mathbf{q}_{d}\mathbf{q}^{T}\}\mathbf{Q}}{\sqrt{2}\sqrt{1{+}\mathbf{q}^{T}\mathbf{q}_{d}}}\\
\boldsymbol{\Xi}&{=}&\frac{\{({}^{b}\mathbf{e}_{q}\mathbf{q}^{T}_{d})(\mathbf{q})^{\times}+({}^{b}\mathbf{e}_{q}\mathbf{q}^{T})(\mathbf{q}_{d})^{\times}\}\mathbf{Q}_{d}}{2(1{+}\mathbf{q}^{T}\mathbf{q}_{d})}
\end{IEEEeqnarray}
Alternative expression for \eqref{eq:d},
\begin{IEEEeqnarray}{rCl}
\mathbf{d}&{=}&-(\mathbf{Q}^{T}\mathbf{Q}_{d}{}^{b}\boldsymbol{\omega}_{d})^{\times}{}^{b}\mathbf{e}_{\omega}-\mathbf{Q}^{T}\mathbf{Q}_{d}{}^{b}\dot{\boldsymbol{\omega}}_{d}\label{eq:dalt}
\end{IEEEeqnarray}

\end{document}